\documentclass[12pt]{amsart}
\usepackage{amscd, amssymb}

\theoremstyle{plain}
\newtheorem{proposition}{Proposition}
\newtheorem{lemma}{Lemma}
\newtheorem{theorem}{Theorem}

\newtheorem{corollary}{Corollary}
\newtheorem*{claim*}{Claim}

\newtheorem*{conjecture*}{Conjecture}
\newtheorem*{question*}{Question}
\newtheorem{remark}{Remark}

\theoremstyle{definition}
\newtheorem{definition}{Definition}
\newtheorem*{definition*}{Definition}

\theoremstyle{remark}
\newtheorem*{remark*}{Remark}
\newtheorem*{remarks*}{Remarks}

\newtheorem*{example*}{Example}
\newtheorem*{acknowledgments}{Acknowledgments}

\newcommand{\PP}{\mathbb{P}}
\newcommand{\QQ}{\mathbb{Q}}
\newcommand{\CC}{\mathbb{C}}
\newcommand{\Mbar}{{\overline{M}}}

\newcommand{\p}{\partial}

\newcommand{\ev}{{\operatorname{ev}}}
\newcommand{\ft}{\operatorname{ft}}

\newcommand{\vir}{\operatorname{vir}}
\newcommand{\la}{\langle}
\newcommand{\ra}{\rangle}
\newcommand{\on}{\operatorname}
\newcommand{\HH}{\mathcal{H}}

\title[Invariance of tautological equations II]
{Invariance of tautological equations II:\\ Gromov--Witten theory}

\author{Y.-P.~Lee}

\address{Department of Mathematics \\
        University of Utah \\
        Salt Lake City, Utah 84112-0090\\
        U.S.A.}
\email{yplee@math.utah.edu}

\thanks{Research partially supported by NSF and AMS Centennial Fellowship}

\begin{document}

\begin{abstract}
The aim of Part II is to explore the technique of \emph{invariance 
of tautological equations} in the realm of Gromov--Witten theory.
The main result is a proof of Invariance Theorem 
(Invariance Conjecture~1 in \cite{ypL1}), 
via the techniques from Gromov--Witten theory.
It establishes some general inductive structure of the tautological rings,
and provides a new tool to the study of this area.
\end{abstract}

\maketitle

\setcounter{section}{-1}

\section{Introduction} \label{s:0}

This work is a continuation of Part I of ITE \cite{ypL1}.
The purpose of this paper is to explore the technique of \emph{invariance 
of tautological equations} in the realm of Gromov--Witten theory.

In Part I, a set of three conjectures on the structure of the tautological
rings were proposed.
The main focus there was to study the \emph{linear} invariance operators
\begin{equation} \label{e:rl}
 \mathfrak{r}_l: R^{k}(\Mbar^{\bullet}_{g,n}) \to 
  R^{k+l-1}(\Mbar^{\bullet}_{g-1,n+2})
\end{equation}
between the tautological rings of moduli of curves,
and their conjectural implications.
Here $\bullet$ stands for possibly disconnected curves.
Note that the arithmetic genus for disconnected curve is defined to be
\[
  g(C) := \sum_{i=1}^d g(C_i) - d +1,
\]
where $C_i$ are connected components of $C$, $C= \amalg_{i=1}^d C_i$.
By the definition given in \cite{ypL1}, the curves in the 
image of $\mathfrak{r}_l$ have at most one more connected components.
Therefore, the connected components of the image curves would have either
smaller genus or the same genus but less marked points 
than those of the domain curves.
Furthermore, Invariance Conjecture~2 asserts that the product of
$\mathfrak{r}_l$ for finitely number of $l$ will be \emph{injective}.
This implies, conjecturally, that the tautological rings have a
previously unknown inductive structure.

The definition of $\mathfrak{r}_l$ given in Part I is via
operations on the decorated graphs.
Note that one class in $R^k(\Mbar_{g,n})$ may have more than one
graphical presentations due to the existence of tautological relations.
It is highly nontrivial that certain combination of these graphical operations
would descend to operations on $R^k(\Mbar_{g,n})$.
The main result (Theorem~\ref{t:c1}) of this paper is to prove the 
the invariance operators $\mathfrak{r}_l$ is well-defined on
$R^k(\Mbar_{g,n})$.  
This will be called \emph{Invariance Theorem}, following \cite{rV}.
The existence alone easily implies some new results in tautological rings
as well as simplified proofs of old ones.
This will be explained in Section~\ref{s:6.5}.

The techniques for the proof comes from interactions of moduli of curves
and Gromov--Witten theory.
The idea of using Gromov--Witten theory to study tautological rings on 
moduli of curves is not new.
The fixed point loci of localizations on moduli of maps
are moduli of curves, and ``trivial'' identities on moduli of maps can
produce non-trivial identities on their fixed point loci.
See \cite{rV} for a nice survey of this subject.
Here however Gromov--Witten theory is used in a different way.
Roughly, instead of localization, 
deformation theory of Gromov--Witten theory is used.
Along the way, some results in Gromov--Witten theory are also proved.

Here is a summary of the content of this paper.
Section~\ref{s:1} gives a quick summary of geometric Gromov--Witten theory.
As alluded in Part I \cite{ypL1}, the set of conjectures proposed there 
are motivated by study of Givental's axiomatic Gromov--Witten theory. 
In Sections~\ref{s:2} and \ref{s:3}, we summarize Givental's theory.
Some geometric results on the tautological classes, known to 
experts, are given in Section~\ref{s:4} for readers' convenience.
In Section~\ref{s:5}, it is proved that each tautological \emph{class} 
is invariant under the quantized lower triangular loop groups.
Thus the invariance under lower triangular loop groups gives no constraints.
In Section~6, we study the invariance of tautological \emph{equations}
under the quantized upper triangular loop groups.
It turns out that this invariance poses very strong constraints on
possible forms of the tautological equations.
As a matter of fact, \emph{the invariance constraints are so strong as to, 
conjecturally, uniquely determine the tautological equations}.
This is the motivation of the three Invariance Conjectures 
advanced in \cite{ypL1}.
The main result, Invariance Theorem (Theorem~\ref{t:c1}), is proved there.

At the end of the paper, we indicate a few applications of Invariance Theorem,
including a Faber type statement for the tautological rings \cite{AL3} and
a uniform derivation of all known tautological equations \cite{GL, AL1, AL2}.
An appendix (jointly with Y.~Iwao) demonstrates some properties of the 
upper triangular loop groups associated to $\PP^1$.

\begin{acknowledgments}
I wish to thank D.~Arcara, E.~Getzler, A.~Givental, R.~Pandharipande,
C.~Teleman and R.~Vakil for many useful discussions.
\end{acknowledgments}

\section{Geometric Gromov--Witten theory} \label{s:1}

\subsection{Preliminaries of Gromov--Witten theory} \label{s:1.1}
Gromov--Witten theory studies the tautological intersection theory on
$\Mbar_{g,n}(X,\beta)$, the moduli spaces of stable maps from curves $C$ of
genus $g$ with $n$ marked points to a smooth projective variety $X$. The
intersection numbers, or \emph{Gromov--Witten invariants}, are integrals of
tautological classes over the virtual fundamental classes of
$\Mbar_{g,n}(X,\beta)$
\[
  \int_{[\Mbar_{g,n}(X,\beta)]^{\vir}}
    \prod_{i=1}^n \ev_i^*(\gamma_i) \psi_i^{k_i}.
\]
Here $\gamma_i \in H^*(X)$
and $\psi_i$ are the \emph{cotangent classes} (gravitational
descendents).

For the sake of the later reference, let us fix some notations.
\begin{enumerate}
\item[(i)] $H := H^*(X, \QQ)$ is a $\QQ$-vector space, assumed of rank $N$.
Let $\{ \phi_{\mu} \}_{\mu=1}^N$ be a basis of $H$.

\item[(ii)] $H$ carries a symmetric bilinear form, Poincar\'e pairing,
\[
  \langle \cdot, \cdot \rangle : H \otimes H \to \QQ.
\]
Define 
\[
  g_{\mu \nu} := \langle \phi_{\mu}, \phi_{\nu} \rangle
\]
and $g^{\mu \nu}$ to be the inverse matrix.

\item[(iii)] Let $\mathcal{H}_t := \oplus_{k=0}^{\infty} H$ be the infinite
dimensional complex vector space with basis $\{ \phi_{\mu} \psi^k \}$.
$\mathcal{H}_t$ has a natural $\QQ$-algebra structure:
\[
 \phi_{\mu} \psi^{k_1} \otimes \phi_{\nu} \psi^{k_2} \mapsto 
 (\phi_{\mu} \cdot \phi_{\nu})  \psi^{k_1 + k_2},
\]
where $\phi_{\mu} \cdot \phi_{nu}$ is the cup product in $H$.

\item[(iv)] Let $\{ t^{\mu}_k \}$, $\mu=1, \ldots, N$, $k=0, \ldots, \infty$,
be the dual coordinates of the basis $\{ \phi_{\mu} \psi^k \}$.

\end{enumerate}

We note that at each marked point, the insertion is $\mathcal{H}_t$-valued. 
Let
\[
 t:= \sum_{k, \mu} t^{\mu}_k \phi_{\mu} \psi^k
\] 
denote a general element in the vector space $\mathcal{H}_t$.
To simplify the notations, $t_k$ will stand for the vector
$(t^1_{k},\ldots,t^N_{k})$ and $t^{\mu}$ for
$(t^{\mu}_0, t^{\mu}_1, \ldots )$.

\begin{enumerate}

\item[(v)] Define 
\[ 
 \la \p^{\mu_1}_{k_1} \ldots \p^{\mu_n}_{k_n}
 \ra_{g,n,\beta} := \int_{[\Mbar_{g,n}(X,\beta)]^{\vir}} \prod_{i=1}^n
 \ev_i^*(\phi_{\mu_i}) \psi_i^{k_i}
\]
and define 
\[
 \la t^n \ra_{g,n,\beta}=\la t \ldots t \ra_{g,n,\beta}
\]
by multi-linearity.

\item[(vi)] Let 
\[
  F^X_g(t) := \sum_{n, \beta} \frac{1}{n!} \la t^n \ra_{g,n,\beta}
\]
be the generating function of all genus $g$ Gromov--Witten invariants.
\footnote{In Gromov--Witten theory, one usually has to deal with 
the coefficients in the \emph{Novikov ring}.
We shall not touch upon this subtleties here but refer the readers to
\cite{LP}.}
The \emph{$\tau$-function of $X$} is the formal expression 
\begin{equation} \label{e:1}
 \tau_{GW}^X := e^{\sum_{g=0}^{\infty} \hbar^{g-1} F_g^X}.
\end{equation}
\end{enumerate}

\subsection{Gravitational ancestors and the $(3g-2)$-jet properties} 
\label{s:1.2}
Let
\begin{equation} \label{e:st}
  \on{st}: \Mbar_{g,m+l}(X,\beta) \to \Mbar_{g,m+l}
\end{equation}
be the \emph{stabilization morphism} defined by forgetting the map and
\[
  \on{ft}: \Mbar_{g,m+l} \to \Mbar_{g,m}
\]
be the \emph{forgetful morphism} defined by forgetting the last $l$ points.
The \emph{gravitational ancestors} are defined to be
\begin{equation} \label{e:ancestor}
  \bar{\psi}_i := (\on{ft} \circ \on{st})^* \psi_i
\end{equation}
and genus $g$ ancestor potential is defined by
\[
  \overline{F}^X_g (t, s) := \sum_{m,l,\beta} \frac{Q^{\beta}}{m! l!}
     \int_{[\Mbar_{g,m+l}(X,\beta)]^{\vir}} 
     \prod_{i=1}^m \sum_{k} t_k^{\mu} (\bar{\psi}_i)^k  \ev_i^*(\phi_{\mu})
     \prod_{i=m+1}^{m+l} \sum_{\mu} s^{\mu} \ev_i^*(\phi_{\mu}) .
\]
The following property is called the $(3g-2)$-jet property \cite{eG4}
\begin{equation} \label{e:3g-2}
  \frac{\p^m}{\p t^{\mu_1}_{k_1 +1} \ldots \p t^{\mu_m}_{k_m +1}}
    \overline{F}^X_g (t, s) |_{t_0 =0} = 0 \qquad 
    \text{for} \ \sum k_i \ge 3g-2.
\end{equation}
This follows from the dimension counting
\[
  \dim \Mbar_{g,n} = 3g-3 + n.
\]

The ancestors and descendents are different, but easy to compare.
Let $D_j$ be the (virtual) divisor on
$\Mbar_{g,n+m}(X,\beta)$ defined by the image of the gluing morphism
\[
 \sum_{\beta'+\beta''=\beta} \sum_{m' + m''=m}
 \Mbar^{(j)}_{0, 2+m'}(X,\beta') \times_X \Mbar_{g,n+m''}(X,\beta'')
 \to \Mbar_{g,n+m}(X,\beta),
\]
where $\Mbar_{g,n+m''}(X,\beta'')$ carries all first $n$ marked points
except the $j$-th one, which is carried by $\Mbar^{(j)}_{0, 2+m'}(X,\beta')$.
It is proved in \cite{KM} that 
\begin{equation} \label{e:compare}
 \psi_j - \bar{\psi}_j = [D_j]^{\vir}.
\end{equation}

\section{Genus zero axiomatic Gromov--Witten theory} \label{s:2}

Let $H$ be a $\QQ$-vector space of dimension $N$ with a distinguished element
$\mathbf{1}$.
Let $\{ \phi_\mu \}$ be a basis of $H$ and $\phi_{\mathbf{1}} = \mathbf{1}$. 
Assume that $H$ is endowed with a nondegenerate symmetric $\QQ$-bilinear form,
or metric, $\langle \cdot,\cdot \rangle$. 
Let $\HH$ denote the infinite dimensional vector space $H[z,z^{-1}]$ 
consisting of Laurent polynomials with coefficients in $H$. 
\footnote{Different completions of $\HH$ are used in different places.
This will be not be discussed details in the present article.
See \cite{LP} for the details.} 
Introduce a symplectic form $\Omega$ on $\HH$:
\[
  \Omega (f(z),g(z) ) := \on{Res}_{z=0} \langle f(-z), g(z) \rangle,
\]
where the symbol $\on{Res}_{z=0}$ means to take the residue at $z=0$.

There is a natural polarization $\HH = \HH_{q}\oplus \HH_{p}$ by the 
Lagrangian subspaces $\HH_{q} := H[z]$ and $\HH_{p} := z^{-1} H [ z^{-1} ]$
which provides a symplectic identification of $(\HH, \Omega)$ with the 
cotangent bundle $T^*\HH_{q}$ with the natural symplectic structure.
$\HH_{q}$ has a basis 
\[
 \{ \phi_{\mu} z^k \}, \quad 1 \le \mu \le N, \quad 0 \le k
\]
with dual coordinates $\{ q^k_\mu \}$.
The corresponding basis for $\HH_{p}$ is
\[
 \{ \phi_{\mu} z^{-k-1} \}, \quad 1 \le \mu \le N, \quad 0 \le k
\]
with dual coordinates $\{ p^k_\mu \}$.

For example, let $\{ \phi_{i} \}$ be an \emph{orthonormal} basis of $H$. 
An $H$-valued Laurent formal series can be written in this basis as
\begin{multline*}
 \ldots + (p^1_1,\ldots,p^N_{1}) \frac{1}{(-z)^2}
 + (p^1_{0},\ldots, p^N_{0}) \frac{1}{(-z)} \\
 + (q^1_{0},\ldots, q^N_{0})
 + (q^1_{1}, \ldots, q^N_{1}) z + \ldots.
\end{multline*}
In fact, $\{ p_k^{i}, q_k^{i} \}$ for $k= 0, 1, 2, \ldots$ and
$i=1, \ldots, N$ are the Darboux coordinates compatible with this
polarization in the sense that
\[
 \Omega = \sum_{i,k} d p^{i}_k \wedge d q^{i}_k .
\]

The parallel between $\HH_q$ and $\HH_t$ is evident, and is in fact
given by the following affine coordinate transformation, 
called the \emph{dilaton shift},
\[
  t^{\mu}_k = q^{\mu}_k + \delta^{\mu \mathbf{1}} \delta_{k 1}.
\]

\begin{definition} \label{d:1} 
Let $G_0(t)$ be a (formal) function on $\HH_t$. 
The pair $T:=(\HH, G_0)$ is called a \emph{$g=0$ axiomatic theory} if 
$G_0$ satisfies three sets of genus zero
tautological equations: 
the {\em Dilaton Equation} \eqref{e:de},
the {\em String Equation} \eqref{e:se} and 
the \emph{Topological Recursion Relations} (TRR) \eqref{e:g0trr}.

\begin{align}
 \label{e:de}
 &\frac{\p G_0(t)}{\p t^{\mathbf{1}}_1} (t)
 = \sum_{k=0}^{\infty} \sum_{\mu} t^{\mu}_k
 \frac{\p G_0(t)}{\p t^{\mu}_k} - 2G_0 (t), \\
  \label{e:se}
  &\frac{\p G_0 (t)}{\p t^{\mathbf{1}}_0} =
  \frac{1}{2} \langle t_0,t_0 \rangle + \sum_{k=0}^{\infty}
  \sum_{\nu} t_{k+1}^{\nu} \frac{\p G_0 (t)}{\p t^{\nu}_k},  \\
  \label{e:g0trr}
  &\frac{\p^3 G_0 (t)}
       {\p t^{\alpha}_{k+1} \p t^{\beta}_{l} \p t^{\gamma}_m}
  = \sum_{\mu \nu} \frac{\p^2 G_0 (t)}{\p t^{\alpha}_{k} \p t^{\mu}_{0}}
    g^{\mu \nu} \frac{\p^3 G_0 (t)}
       {\p t^{\nu}_{0} \p t^{\beta}_{l} \p t^{\gamma}_m}, 
  \quad \forall \alpha, \beta, \gamma, k,l,m.
\end{align}
\end{definition}

In the case of geometric theory, $G_0 = F_0^X$ 
It is well known that $F_0^X$ satisfies the above three
sets of equations \eqref{e:de} \eqref{e:se} \eqref{e:g0trr}.
The main advantage of viewing the genus zero theory through this formulation, 
seems to us, is to replace $\mathcal{H}_t$ by $\HH$ where a symplectic 
structure is available. 
Therefore many properties can be reformulated in terms of the symplectic
structure $\Omega$ and hence independent of the choice of the polarization.
This suggests that the space of genus zero axiomatic Gromov--Witten
theories, i.e. the space of functions $G_0$ satisfying the string equation, 
dilaton equation and TRRs, has a huge symmetry group.

\begin{definition} \label{d:2}
Let $L^{(2)}GL(H)$ denote the {\em twisted loop group} which consists of
$\operatorname{End}(H)$-valued formal Laurent series $M(z)$ in
the indeterminate $z^{-1}$ satisfying $M^*(-z)M(z)=\mathbf{I}$.
Here $\ ^*$ denotes the adjoint with respect to $(\cdot ,\cdot )$.
\end{definition}

The condition $M^*(-z)M(z)=\mathbf{I}$ means that $M(z)$ is a 
symplectic transformation on $\HH$.

\begin{theorem} \label{t:1} \cite{aG4}
The twisted loop group acts on the space of axiomatic genus zero theories.
Furthermore, the action is transitive on the semisimple theories of a fixed
rank $N$.
\end{theorem}

\begin{remarks*}
(i) In the geometric theory, $F^X_0(t)$ is usually a formal function in $t$.
Therefore, the corresponding function in $q$ would be formal at 
$q = - \mathbf{1}z$.
Furthermore, the Novikov rings are usually needed to ensure the 
well-definedness of $F^X_0(t)$. (cf.~Footnote~1.)

(ii) It can be shown that the axiomatic genus zero theory over complex numbers
is equivalent to the definition of abstract (formal) Frobenius manifolds, 
not necessarily conformal. The coordinates on the
corresponding Frobenius manifold is given by the following map \cite{DW}
\begin{equation} \label{e:dg}
 s^{\mu} := \frac{\p}{\p t^{\mu}_0}\frac{\p}{\p t^{\mathbf{1}}_0} G_0 (t).
\end{equation}
\emph{From now on, the term ``genus zero axiomatic theory'' is identified 
with ``Frobenius manifold''}.

(iii) The above formulation (or the Frobenius manifold 
formulation) does not include the divisor axiom, which is true for any
geometric theory.

(iv) Coates and Givental \cite{CG} (see also \cite{aG4}) give a beautiful 
geometric reformation of the genus zero axiomatic theory in terms of 
\emph{Lagrangian cones} in $\HH$.
When viewed in the Lagrangian cone formulation, Theorem~\ref{t:1} 
becomes transparent and a proof is almost immediate.
\end{remarks*}

\section{Quantization and higher genus axiomatic theory} \label{s:3}

\subsection{Preliminaries on quantization}

To quantize an infinitesimal symplectic transformation, or its corresponding
quadratic hamiltonians, we recall the standard Weyl quantization. 
A polarization $\HH=T^* \HH_q$ on the symplectic vector space $\HH$ (the phase
space) defines a configuration space $\HH_q$. 
The quantum ``Fock space'' will be a certain class of functions 
$f(\hbar, q)$ on $\HH_q$ (containing at least polynomial functions), 
with additional formal variable $\hbar$ (``Planck's constant''). 
The classical observables are certain functions of $p, q$. The
quantization process is to find for the classical mechanical system on $\HH$ a
``quantum mechanical'' system on the Fock space such that the classical
observables, like the hamiltonians $h(q,p)$ on $\HH$, are quantized to become
operators $\widehat{h}(q,\dfrac{\p}{\p q})$ on the Fock space.

Let $A(z)$ be an $\on{End}(H)$-valued Laurent formal series in $z$ satisfying
\[
  (A(-z) f(-z), g(z)) + (f(-z), A(z) g(z)) =0,
\]
then $A(z)$ defines an infinitesimal symplectic transformation
\[
  \Omega(A f, g) + \Omega(f, A g)=0.
\]
An infinitesimal symplectic transformation $A$ of $\HH$ corresponds to a
quadratic polynomial $P(A)$ in $p, q$
\[
  P(A)(f) := \frac{1}{2} \Omega(Af, f) .
\]

Choose a Darboux coordinate system $\{ q^i_k, p^i_k \}$. 
The quantization $P \mapsto \widehat{P}$ assigns
\begin{equation} \label{e:wq}
 \begin{split}
  &\widehat{1}= 1, \  \widehat{p_k^{i}}= \sqrt{\hbar} \frac{\p}{\p q^{i}_k}, \
   \widehat{q^{i}_k} = q^{i}_k / {\sqrt{\hbar}}, \\
  &\widehat{p^{i}_k p^{j}_l} = \widehat{p^{i}_k} \widehat{p^{j}_l}
    =\hbar \frac{\p}{\p q^{i}_k} \frac{\p}{\p q^{j}_l}, \\
   &\widehat{p^{i}_k q^{j}_l} = q^{j}_l \frac{\p}{\p q^{i}_k},\\
  &\widehat{q^{i}_k q^{j}_l} = {q}^{i}_k {q}^{j}_l /\hbar ,
 \end{split}
\end{equation}
In summary, the quantization is the process
\[
\begin{matrix}
  A   &\mapsto & P(A)  &\mapsto & \widehat{P(A)} \\
  \text{inf. sympl. transf.}  &\mapsto & \text{quadr. hamilt.}  
    &\mapsto & \text{operator on Fock sp.}.
\end{matrix}
\]
It can be readily checked that the first map is a Lie algebra isomorphism:
The Lie bracket on the left is defined by $[A_1, A_2]=A_1 A_2 - A_2 A_1$ and 
the Lie bracket in the middle is defined by Poisson bracket
\[
 \{ P_1(p,q), P_2(p,q) \} = \sum_{k,i} 
 \frac{\p P_1}{\p p^{i}_k} \frac{\p P_2}{\p q^{i}_k}
 -\frac{\p P_2}{\p p^{i}_k} \frac{\p P_1}{\p q^{i}_k}.
\]
The second map is not a Lie algebra homomorphism, but is very close to being
one.
\begin{lemma} \label{l:1}
\[
 [\widehat{P_1},\widehat{P_2}] =  \widehat{\{ P_1, P_2 \}} + 
  \mathcal{C}(P_1,P_2),
\]
where the cocycle $\mathcal{C}$, in orthonormal coordinates, vanishes except
\[
 \mathcal{C}(p^{i}_k p^{j}_l, q^{i}_k q^{j}_l) = 
 -\mathcal{C}(q^{i}_k q^{j}_l, p^{i}_k p^{j}_l) 
 = 1 + \delta^{i j} \delta_{kl}.
\]
\end{lemma}

\begin{example*}
Let $\dim H=1$ and $A(z)$ be multiplication by $z^{-1}$. 
It is easy to see that $A(z)$ is infinitesimally symplectic.
\begin{equation} \label{e:cse}
 \begin{split}
 P(z^{-1})= &-\frac{q_0^2}{2} - \sum_{m=0}^{\infty} q_{m+1} p_m \\
 \widehat{P(z^{-1})} = &-\frac{q_0^2}{2} 
            - \sum_{m=0}^{\infty} q_{m+1} \frac{\p}{\p q_m}.
 \end{split}
\end{equation}
\end{example*}

Note that one often has to quantize the symplectic instead of the 
infinitesimal symplectic transformations. 
Following the common practice in physics, define
\begin{equation} \label{e:q}
  \widehat{e^{A(z)}} := e^{\widehat{A(z)}},
\end{equation}
for $e^{A(z)}$ an element in the twisted loop group.

\subsection{$\tau$-function for the axiomatic theory}

Let $X$ be the space of $N$ points and $H^{N pt}:= H^*(X)$.
Let $\phi_{i}$ be the delta-function at the $i$-th point.
Then $\{ \phi_{i} \}_{i=1}^N$ form an orthonormal basis 
and are the idempotents of the quantum product
\[
  \phi_{i} * \phi_{j} = \delta_{ij} \phi_{i}.
\]
The genus zero potential for $N$ points is nothing but a sum of genus zero 
potentials of a point
\[
  F^{N pt}_0 (t^1, \ldots, t^N) = F^{pt}_0 (t^1) + \ldots + 
  F^{pt}_0(t^N).
\]
In particular, the genus zero theory of $N$ points is semisimple.

By Theorem~\ref{t:1}, any \emph{semisimple} genus zero axiomatic theory $T$ 
of rank $N$ can be obtained from $H^{N pt}$ by action of an an element $O^T$ 
in the twisted loop group. 
By Birkhoff factorization, $O^T = S^T (z^{-1}) R^T (z)$, where
$S(z^{-1})$ (resp.~$R(z)$) is an matrix-valued functions in $z^{-1}$
(resp.~$z$). 

In order to define the axiomatic higher genus potentials $G_g^T$ for the 
semisimple theory $T$, one first introduces the ``$\tau$-function of $T$''.

\begin{definition} \label{d:3} \cite{aG2}
Define the \emph{axiomatic $\tau$-function} as
\begin{equation} \label{e:taug}
  \tau_G^T:= \widehat{S^T} (\widehat{R^T} \tau^{N pt}_{GW}),
\end{equation}
where $\tau^{N pt}_{GW}$ is defined in \eqref{e:1}.
Define the \emph{axiomatic genus $g$ potential} $G_g^T$ via the formula 
(cf.~\eqref{e:1})
\begin{equation} \label{e:2}
  \tau_G^T =: e^{\sum_{g=0}^{\infty} \hbar^{g-1} G_g^T}.
\end{equation}
\end{definition}

\begin{remark*}
(i) It is not obvious that the above definitions make sense.
The function $\widehat{S^T} (\widehat{R^T} \tau^{N pt})$ is well-defined, 
due to the $(3g-2)$-jet properties \eqref{e:3g-2}, 
proved in \cite{eG4} for geometric Gromov--Witten theory 
and in \cite{aG2} in the axiomatic framework.
The fact $\log \tau_G^T$ can be written as 
$\sum_{g=0}^{\infty} \hbar^{g-1} (\text{formal function in $t$})$ is also
nontrivial.
The interested readers are referred to the original article \cite{aG2}
or \cite{LP} for details.
\end{remark*}

An immediate question regarding Definition~\ref{d:3}:
When the axiomatic semisimple theory actually comes from a projective variety 
$X$, is the $\tau$-function defined in the axiomatic theory
the same as the $\tau$-function defined in the geometric Gromov--Witten theory?
This is known as Givental's conjecture.
Givental himself establishes the special cases when $X$ is toric Fano.
Recently, C.~Teleman \cite{cT} has announced a complete proof based on
his classification of semisimple $2D$ cohomological field theories.
\footnote{Teleman also informed us that M.~Kontsevich has a related
result, which gives a full description of the deformation theory of
(semisimple) $2D$ cohomological field theory.}

\begin{theorem} \label{t:cT} (\cite{aG2} \cite{cT})
Let $X$ be a projective variety whose quantum cohomology is semisimple,
then
\begin{equation} \label{e:cT}
  \tau_G^X = \tau_{GW}^X.
\end{equation}
\end{theorem}

\section{Tautological equations in Gromov--Witten theory} \label{s:4}

As stated in the Introduction, the material here is known to experts,
and is included for lack of a good general reference (known to us).
Therefore, the discussions will be brief.

Due to the stabilization morphism \eqref{e:st}, 
any \emph{tautological equation} in $\Mbar_{g,n}$ can be pulled back and
become an equation on tautological classes on $\Mbar_{g,n}(X,\beta)$.
The $\psi$-classes can be either transformed to ancestor classes $\bar{\psi}$
or to $\psi$-classes on the moduli spaces of stable maps.
Due to the functorial properties of the virtual fundamental classes, 
the pull-backs of the tautological equations hold for Gromov--Witten
theory of any target space.
The term \emph{tautological equations} will also be used for the corresponding
equations in Gromov--Witten theory and in the theory of spin curves.
Equations in Gromov--Witten theory which are valid for all target spaces
are called the \emph{universal equations}.
Tautological equations are universal equations.

However, these induced equations will produce relations among generalized 
Gromov--Witten invariants, which involves not only $\psi$-classes but also
$\kappa$-classes and boundary classes.
Since the integration over boundary classes can be written in terms of
ordinary Gromov--Witten invariants by the splitting axiom, what is really
in question is the $\kappa$-classes. 
We will start with some results in tautological classes on moduli of curves.

Let 
\[
 \ft_l: \Mbar_{g,n+l} \to \Mbar_{g,n}
\]
be the forgetful morphism, forgetting the last $l$ marked points.
\begin{lemma} \label{l:2}
\[
 (\ft_l)_* ( \prod_{i=1}^{n} (\psi_i)^{k_i} 
  \prod_{i+n+1}^{n+l} (\psi_i)^{k_i +1} ) 
 = K_{k_{n+1} \ldots k_{n+l}} \prod_{i=1}^{n} \psi_i^{k_i},
\]
where $K_{k_{n+1} \ldots k_{n+l}}$ is defined as follows.
Let $\sigma \in S_l$ be an element in the symmetric group permuting
the set $\{ n+1, \ldots n+l \}$.
$\sigma$ can be written as a product of disjoint cycles 
\[
 \sigma = c^1 c^2 \ldots, \quad \text{where} \quad
 c^i:= (c^i_1 \ldots c^i_{a^i})
\] 
with $c^i_j \in \{ n+1, \ldots n+l \}$.
Define
\[
  K_{k_{n+1} \ldots k_{n+l}} := 
  \sum_{\sigma \in S_l} \prod_{i (\sigma = \prod c^i)}  
  \sum_{j=1}^{a^i} \kappa_{c^i_j}.
\]
\end{lemma}

For example, when $l=2$, the formula becomes
\[
 \begin{split}
  &(\ft_2)_*  \left( (\prod_{i=1}^{n} (\psi_i)^{k_i}) 
   (\psi_{n+1})^{k_{n+1} +1} (\psi_{n+2})^{k_{n+2} +1} \right) \\
 = &(\kappa_{k_{n+1}} \kappa_{k_{n+2}} + \kappa_{k_{n+1}+k_{n+2}})
   \prod_{i=1}^{n} (\psi_i)^{k_i}.
 \end{split}
\]

\begin{proof}
The proof follows from induction on $l$ and the following three geometric
ingredients, which are well-known in the theory of moduli of curves.
\begin{itemize}
\item Let $\ft_1:\Mbar_{g,n+1} \to \Mbar_{g,n}$ be the forgetful morphism,
forgetting the last marked point and $D_{i,n+1}$ be the boundary divisor
in $\Mbar_{g,n+1}$ defined as the image of the section of $\ft_1$, 
considered as the universal curve, by the $i$-th marked point. Then
\begin{equation} \label{e:18}
  \ft_1^*(\psi_i) = \psi_i - D_{i,n+1}.
\end{equation}

\item $\psi_i$ for $i \le n$ on $\Mbar_{g,n+1}$ vanishes when restricted to
$D_{i,n+1}$.

\item 
\begin{equation} \label{e:19}
  \ft_1^* (\kappa_l) = \kappa_l - \psi_{n+1}^l.
\end{equation}
\end{itemize}
\end{proof}

The following result follows by combining the above ingredients and 
induction on power of $\kappa$-classes.

\begin{corollary} \label{c:1}
A tautological equations on $\Mbar_{g,n}$ involving $\kappa$-classes of 
highest power less than $l$ (e.g.~$\kappa_{k_1} \ldots \kappa_{k_l}$) 
can be written as a pushforward, via forgetful morphism $\ft^l$, 
of a tautological equation on $\Mbar_{g,n+l}$, involving only boundary
strata and $\psi$-classes.
\end{corollary}

By pulling-back to moduli of stable maps, one has the following corollary.

\begin{corollary}
(i) The \emph{system} of generalized Gromov--Witten invariants involving 
$\kappa$-classes and $\lambda$-classes
is the same as the system of usual Gromov--Witten invariants.

(ii) Any induced equation of generalized Gromov--Witten invariants can be 
written as an equation of ordinary Gromov--Witten invariants.
\end{corollary}

\begin{proof}
The fact that the system of generalized GW invariants involving 
$\lambda$-classes can be reduced to ordinary GW invariants follows from 
\cite{FP}. 
The part involving $\kappa$-classes follows from Lemma~\ref{l:2}.
\end{proof}

With this Corollary, one can talk about the induced equations of 
(ordinary) Gromov--Witten invariants from any tautological equations.

\begin{remark} \label{r:compatible}
It is not difficult to see, from the above discussions 
(in particular equations \eqref{e:18}, \eqref{e:19} and Lemma~\ref{l:2}), 
that the three graphical operations introduced
in \cite{ypL1} (cutting edges, genus reduction, and splitting vertices) 
are compatible with the pull-back operations.
\end{remark}

\section{Invariance under lower triangular subgroups} \label{s:5}
The twisted loop group is
generated by ``lower triangular subgroup'' and the ``upper triangular
subgroup''. The lower triangular subgroup consists of $\on{End}(H)$-valued
formal series $S(z^{-1})= e^{s(z^{-1})}$ in $z^{-1}$ satisfying $S^*(-z) S(z) =
\mathbf{1}$ or equivalently
\[
 s^*(-z^{-1}) + s(z^{-1}) =0.
\]

\subsection{Quantization of lower triangular subgroups}

The quadratic hamiltonian of $s(z^{-1}) = \sum_{l=1}^{\infty} s_l z^{-l}$ is
\[
  \sum_{l=1}^{\infty} \sum_{n=0}^{\infty} \sum_{i,j}
  (s_{l})_{ij} q^j_{l+n} p^i_n +
  \sum \frac{1}{2} (-1)^n (s_{l})_{ij} q^i_n q^j_{l-n-1}.
\]
The fact that $s(z^{-1})$ is a series in $z^{-1}$ implies that the
quadratic hamiltonian $P(s)$ of $s$ is of
the form $q^2\text{-term} + qp\text{-term}$ where $q$ in $qp$-term
does not contain $q_0$. The quantization of the $P(s)$
\[
 \hat{s} = \sum (s_{l})_{ij} q^j_{l+n} \p_{q^i_n} +
 \frac{1}{2\hbar} \sum (-1)^{n} (s_{l})_{ij} q^i_n q^j_{l-n-1} .
\]
Here $i,j$ are the indices of the orthonormal basis. (The indices $\mu, \nu$
will be reserved for the ``gluing indices'' at the nodes.)
For simplicity of the notation, we adopt the \emph{summation convention} to
\emph{sum over all repeated indices}.

Let {$\displaystyle \frac{d \tau_G}{d \epsilon_s} := \hat{s}(z) \tau_G$}. Then
\[
 \begin{split}
  \frac{d G_0(\epsilon_s)}{d \epsilon_s} =
   &\sum_{l=1}^{\infty} \sum_{n=0}^{\infty}
        \sum_{i,j} (s_{l})_{ij} q^j_{l+n} \p_{q^i_n} G_0
   + \frac{1}{2} (-1)^n (s_{l+n+1})_{ij} q^i_n q^j_l . \\
  \frac{d G_g(\epsilon_s)}{d \epsilon_s} =
   &\sum_{l=1}^{\infty} \sum_{n=0}^{\infty}
        \sum_{i,j} (s_{l})_{ij} q^j_{l+n} \p_{q^i_n} G_g,
    \quad \text{for $g \ge 1$}. \\
 \end{split}
\]

Define
\[
  \la \p^{i_1}_{k_1} \p^{i_2}_{k_2} \ldots \p^{i_n}_{k_n} \ra_g
  := \frac{\p^n G_g}{\p t^{i_1}_{k_1} \p t^{i_2}_{k_2} \ldots
    \p t^{i_n}_{k_n}},
\]
and denote $\la \ldots \ra := \la \ldots \ra_0$. These functions $\la \ldots
\ra_g$ will be called \emph{axiomatic Gromov--Witten invariants}. Then
\begin{equation} \label{e:s0}
 \begin{split}
 &\frac{d}{d \epsilon_s} \langle \p^{i_1}_{k_1} \p^{i_2}_{k_2} \ldots \rangle\\
 = &\sum (s_l)_{ij} q^j_{l+n} \langle \p^i_n \p^{i_1}_{k_1} \ldots \rangle
  +\sum_{l=1}^{\infty} \sum_{i,a} (s_l)_{i i_a} \langle \p^i_{k_a-l}
   \p^{i_1}_{k_1} \ldots \hat{\p^{i_a}_{k_a}} \ldots \rangle \\
  + &\frac{\delta}{2} \Big( (-1)^{k_1} \sum (s_{k_1+k_2+1})_{i_1 i_2}
   + (-1)^{k_2} \sum (s_{k_1+k_2+1})_{i_2 i_1} \Big),
 \end{split}
\end{equation}
where $\delta=0$ when there are more than 2 insertions and
$\delta=1$ when there are two insertions. The notation
$\hat{\p^{i}_{k}}$ means that $\p^{i}_{k}$ is omitted from the summation.
We \emph{assume that there are at least two insertions},
as this is the case in our application.

For $g\ge 1$
\begin{equation} \label{e:sg}
 \begin{split}
  &\frac{d}{d\epsilon_s} \langle \p^{i_1}_{k_1} \p^{i_2}_{k_2} \ldots \rangle_g\\
  = &\sum (s_l)_{ij} q^j_{l+n} \langle \p^i_n \p^{i_1}_{k_1} \ldots \rangle_g
  + \sum \sum_{a} (s_l)_{i i_a} \langle \p^i_{k_a-l}
   \p^{i_1}_{k_1} \ldots \hat{\p^{i_a}_{k_a}} \ldots \rangle_g
 \end{split}
\end{equation}

\subsection{$S$-Invariance}

\begin{theorem} \label{t:sinv} \emph{($S$-invariance theorem)}
All tautological equations are invariant under action of lower triangular
subgroups of the twisted loop groups.
\end{theorem}

\begin{proof}
Let $E=0$ be a tautological equation of axiomatic Gromov--Witten invariants.
Suppose that this equation holds for a given semisimple Frobenius manifold,
e.g.~$H^{N pt} \cong \CC^N$. We will show that $\hat{s} E =0$. This will prove
the theorem.

$\hat{s} E =0$ follows from the following facts:
\begin{enumerate}
\item[(a)] The combined effect of the first term in \eqref{e:s0} (for genus
zero invariants) and in \eqref{e:sg} (for $g \ge 1$ invariants) vanishes.
\item[(b)] The combined effect of the remaining terms in \eqref{e:s0} and in
\eqref{e:sg} also vanishes.
\end{enumerate}
(a) is due to the fact that the sum of the contributions from the first
term is a derivative of the original equation $E=0$ with respect to $q$
variables. Therefore it vanishes.

It takes a little more work to show (b).
Recall that all tautological equations are induced from moduli spaces of
curves. Therefore, any relations of tautological classes on $\Mbar_{g,n}$
contain no genus zero components of two or less marked points.
However, when one writes down the induced equation for
(axiomatic) Gromov--Witten invariants, the genus zero invariants with two
insertions will appear. This is due to the difference between the cotangent
classes on $\Mbar_{g,n+m}(X, \beta)$ and the pull-back classes from
$\Mbar_{g,n}$. Therefore the only contribution from the third term of
\eqref{e:s0} comes from these terms. More precisely, let $\psi_j$ (descendents)
denote the $j$-th cotangent class on $\Mbar_{g,n+m}(X,\beta)$ and
$\bar{\psi}_j$ (ancestors) the pull-backs of cotangent classes from
$\Mbar_{g,n}$ by the combination of the stabilization and forgetful morphisms
(forgetting the maps and extra marked points, and stabilizing if necessary).

Denote $\la \p^{\mu}_{k,\bar{l}}, \ldots \ra$
the generalized (axiomatic) Gromov--Witten invariants with
$\psi_1^k \bar{\psi}_1^l \ev_1^*(\phi_{\mu})$ at the first marked point.
The equation \eqref{e:compare} can be rephrased in terms of invariants as
\begin{equation*} 
  \la \p^i_{k,\bar{l}} \ldots \ra_g =
  \la \p^i_{k+1, \overline{l-1}} \ldots \ra_g
  - \la \p^i_k \p^{\mu} \ra g_{\mu \mu'} 
  \la \p^{\mu'}_{\overline{l-1}} \ldots \ra_g.
\end{equation*}
For simplicity, denote
\[
  \la \ldots \p^{\mu} \ra \la \p^{\mu} \ldots \ra 
  := \la \ldots \p^{\mu} \ra g_{\mu \mu'} \la \p^{\mu'} \ldots \ra .
\]
Repeat this process of reducing $\bar{l}$, one can show by induction that
\[
 \begin{split}
  &\la \p^i_{k,\bar{l}} \ldots \ra_g =
  \la \p^i_{k+r, \overline{l-r}} \ldots \ra_g
  -\la \p^i_{k+r-1} \p^{\mu_1} \ra
   \la \p^{\mu_1}_{\overline{l-r}} \ldots\ra_g - \ldots \\
  &-\la \p^i_{k} \p^{\mu_1} \ra \left[ \sum_{p=1}^r (-1)^{p+1}
  \sum_{k_1 + \ldots + k_p = r-p} \la \p^{\mu_1}_{k_1} \p^{\mu_2} \ra
  \ldots \la \p^{\mu_{p-1}}_{k_{p-1}} \p^{\mu_p} \ra
  \la \p^{\mu_p}_{k_p, \overline{l-r}} \ldots \ra_g \right] .
 \end{split}
\]

Now suppose that one has an equation of tautological classes of $\Mbar_{g,n}$.
Use the above equation (for $r=l$) one can translate the equation of
tautological classes on $\Mbar_{g,n}$ into an equation of the (axiomatic)
Gromov--Witten invariants. The term-wise cancellation of the contributions
from the second and the third terms of \eqref{e:s0} and \eqref{e:sg} can
be seen easily by straightforward computation.
\end{proof}

If the above description is a bit abstract, the reader might want to
try the following simple example.
$\psi_1^2$ on $\Mbar_{g,1}$ is translated into invariants:
\[
 \la \p^x_2 \ra_g - \la \p_1^x \p^{\mu} \ra \la \p^{\mu} \ra_g
  - \la \p^x \p^{\mu} \ra \la \p_1^{\mu} \ra_g
  + \la \p^x \p^{\mu} \ra \la \p^{\mu} \p^{\nu} \ra \la \p^{\nu} \ra_g .
\]
The above ``translation'' from tautological classes to Gromov--Witten
invariants are worked out explicitly in some examples in Sections~6 and 7 of
\cite{eG2}.

\begin{remark} \label{r:1}
The $S$-invariance theorem actually hold at the level of (Chow or
cohomology) classes, rather than just the numerical invariants.
The geometric content is \eqref{e:compare}.
This should be clear from the proof.

\end{remark}

\subsection{Reduction to $q_0=0$}
The arguments in this section are mostly taken from \cite{GL}.

Let $E=0$ be a tautological equation of (axiomatic) Gromov--Witten invariants.
Since we have already proved $\hat{s}(E)=0$, our next goal would be
to show $\hat{r}(E)=0$. In this section, we will show that it suffices to
check $\hat{r}(E)=0$ on the subspace $q_0=0$.

\begin{lemma}
It suffices to show $\hat{r} E=0$ on each level set of the map $q \mapsto s$
in \eqref{e:dg}.
\end{lemma}

\begin{proof}
The union of the level sets is equal to $\HH_+$.
\end{proof}

\begin{lemma} \label{l:3}
It suffices to check the relation for all $\hat{r}(z) E =0$ along
$z\HH_{+}$ (i.e.~$q_0=0$).
\end{lemma}

\begin{proof}
Theorem~5.1 of \cite{aG2} states that a particular lower triangular 
matrix $S_s$, which is called ``calibration'' of the Frobenius manifold, 
transforms the level set at $s$ to $z \HH_+$. $S$-invariance Theorem then 
concludes the proof.
\end{proof}

\begin{remark} \label{r:2}
In fact, $S_s$ can be taken as a fundamental solution of the horizontal 
sections of the Dubrovin (flat) connection, in $z^{-1}$ formal series. 
It was discovered in \cite{aG2}, 
following the works in \cite{KM} and \cite{eG4}, that 
\[
  \mathcal{A}:= \hat{S}_s \tau^X
\]
is the corresponding generating function for ``ancestors''.
Therefore the transformed equation $\hat{S}_s E \hat{S}_s^{-1} =0$ is really 
an equation of \emph{ancestors}.
\end{remark}

\section{Invariance under upper triangular subgroups} \label{s:6}

\subsection{Quantization of upper triangular subgroups}

The upper triangular subgroup consists of the regular part of the twisted loop
groups $R(z)= e^{r(z)}$ satisfying 
\begin{equation} \label{e:R}
R^*(-z) R(z) = \mathbf{1} 
\end{equation}
or equivalently
\begin{equation} \label{e:r}
  r^*(-z) + r(z) =0.
\end{equation}

The quantization of $r(z)$ is
\[
 \begin{split}
 \hat{r}(z) = &\sum_{l=1}^{\infty} \sum_{n=0}^{\infty}
    \sum_{i,j} (r_l)_{ij} q^j_n \p_{q^i_{n+l}} \\
  + &\frac{\hbar}{2} \sum_{l=1}^{\infty} \sum_{m=0}^{l-1}
   (-1)^{m+1} \sum_{i j} (r_l)_{ij} \p_{q^i_{l-1-m}} \p_{q^j_m}.
 \end{split}
\]

Therefore
\begin{equation} \label{e:rg}
 \begin{split}
  &\frac{d \langle \p^{i_1}_{k_1} \p^{i_2}_{k_2} \ldots \rangle_g}
        {d \epsilon_r} \\
  = &\sum_{l=1}^{\infty} \sum_{n=0}^{\infty} \sum_{i,j} (r_l)_{ij} q^j_n
   \langle \p^i_{n+l} \p^{i_1}_{k_1} \ldots \rangle_g \\
  + &\sum_{l=1}^{\infty} \sum_{i,a} (r_l)_{i i_a} \langle \p^i_{k_a+l}
    \p^{i_1}_{k_1} \ldots \hat{\p^{i_a}_{k_a}} \ldots \rangle_g \\
  + &\frac{1}{2} \sum_{l=1}^{\infty} \sum_{m=0}^{l-1} (-1)^{m+1}
   \sum_{i j} (r_l)_{ij} \langle \p^i_{l-1-m} \p^j_m
    \p^{i_1}_{k_1} \p^{i_2}_{k_2} \ldots \rangle_{g-1} \\
  + &\frac{1}{2} \sum_{l=1}^{\infty} \sum_{m=0}^{l-1} (-1)^{m+1}
   \sum_{i j} \sum_{g'=0}^g (r_l)_{ij} \p^{i_1}_{k_1} \p^{i_2}_{k_2} \ldots
   ( \langle \p^i_{l-1-m} \rangle_{g'} \langle \p^j_m \rangle_{g-g'} ).
 \end{split}
\end{equation}
Here, if $g=0$, the third term on the right $\langle \ldots \rangle_{-1} = 0$
by definition.
Also, it is understood that the formula for $\hat{r}_l$
extends to products of Gromov--Witten invariants by Leibniz rule.

\subsection{Relations to invariance of tautological equations}
Let $E=0$ be a tautological equation on moduli of curves.
As explained in Part I, it can be written in terms of a formal sum of
decorated graphs.
Denote $E=0$ also the induced equation of Gromov--Witten invariants. 
Consider $\frac{d}{d \epsilon_r} E$.
It is clear that the first term of \eqref{e:rg} vanish as
$E=0$ implies $\sum (r_l)_{ij} q^j_n E =0$.
Similarly, the contribution to the second term from $\hat{\p^i_k}$ at an
external marked point (i.e.~not at a node) cancels.
Therefore, $\frac{d}{d \epsilon_r} E$ consists of three parts, 
from the second term (with $\hat{\p^i_k}$ at a node), 
third term and fourth term.

It follows from the usual correspondence between tautological classes and
Gromov--Witten invariants that these three parts corresponds to 
three graphical operations defined in \cite{ypL1}: 
\begin{itemize}
\item The second term, when $i_a$ are indices at a node,
corresponds to \emph{cutting the edges}.
\item The third term corresponds to \emph{genus reduction}.
\item The last term corresponds to \emph{splitting the vertices}.
\end{itemize}
Equation \eqref{e:r} implies $(r_l)_{ij}$ is symmetric in $i,j$ for $l$ even 
and anti-symmetric for $l$ odd.
The corresponding operation for a fixed $l$ on the decorated graphs 
are denoted $\mathfrak{r}_l$ in Part I.
In fact, this is the original motivation of the Invariance Conjecture~1 and 2.
Invariance Conjecture~3 is based upon a folklore belief:
Call an equation of Gromov--Witten invariants ``universal'' if it holds
for all Gromov--Witten theory.
It is believed that all universal equations are induced 
from tautological equations on moduli of curves.
Invariance Conjecture~3 is statement of this speculation in terms of 
invariance constraints.

To show that $\mathfrak{r}_l$ is well-defined is equivalent to showing
\begin{equation} \label{e:rl2}
  \mathfrak{r}_l (E) =0
\end{equation}
for any tautological equation 
\[
 E=0
\]
in $R^k(\Mbar_{g,n})$.
Using the link between $\mathfrak{r}_l$ on moduli of curves and
$\widehat{r}_l$ on Gromov--Witten theory, one's first step is to show that
\begin{equation} \label{e:Ri}
  \widehat{r}_l (\tilde{E}) =0
\end{equation}
for Gromov--Witten theory.
This will serve as a numerical invariant form of Invariance Theorem,
which asserts the invariance at the level of cycles.

Note that equation \eqref{e:Ri} can be interpreted as an infinitesimal form
of the requirement that the induced tautological equation 
\[
 \tilde{E} =0
\]
has the same form for any Gromov--Witten theory.
That is, the tautological equation is invariant under the action of 
quantized upper triangular subgroups of the twisted loop groups.
\emph{This is the reason the term ``invariance of tautological equations''
is used for equation \eqref{e:rl2}.}

\begin{theorem} \label{t:rinv} 
Equation \eqref{e:Ri} holds. 
That is, all tautological equations of Gromov--Witten invariants 
are invariant under the action of upper triangular
subgroups of the twisted loop groups.
\end{theorem}

\begin{proof} 
Since we have already established $S$-invariance theorem,
this theorem follows from Theorem~\ref{t:cT}.
Indeed, \eqref{e:cT} implies that there is a ``loop group'' element 
(or rather its quantization) taking a tautological equation on moduli
of curves to that of any semisimple theory \emph{and vice versa}.
Therefore the set of induced tautological equations of one semisimple theory 
has a one-one correspondence with 
the set of tautological equations of another semisimple theory.
\end{proof}

\subsection{Proof of Invariance Theorem}
A set of three Invariance Conjectures are advanced in \cite{ypL1} to
give the tautological relations an inductive structure.
The main purpose of this section is to establish Conjecture~1 there, 
which corresponds to an infinitesimal form of Theorem~\ref{t:rinv}
on moduli of curves.

\begin{theorem} \label{t:c1}  \emph{(Invariance Theorem)}
$\mathfrak{r}_l$ in \eqref{e:rl} is well-defined, or equivalently
\eqref{e:rl2} holds for any tautological equation.
That is, Invariance Conjecture~1 in \cite{ypL1} is true. 
\end{theorem}

\begin{proof}
With all the preparation above, there are two remaining ingredients 
in the proof:
\begin{enumerate}
\item Teleman's classification theorem of semisimple
cohomological field theories \cite{cT}, 
applied to Givental's framework on $X=\PP^1$.
\item An explicit calculation of $r_l$-matrix of $\PP^1$ in the Appendix.
\end{enumerate}
Note that (1) is stronger than Theorem~\ref{t:cT}:
It implies a ``cycle form'' of Givental's formula.

Let 
\[
  E= \sum_i c_i \Gamma_i = 0 
\]
be a tautological relation. (Notations as in Part I \cite{ypL1}.)
The induced tautological equation on Gromov--Witten theory of two points
is denoted $\tilde{E}$.
(1) implies that this tautological equation on the (two copies of) 
moduli space of curves ($X = 2 pt$) is transformed to the corresponding 
tautological equation on moduli space of stable maps to $\PP^1$.
By Theorem~\ref{t:sinv} and Remark~\ref{r:1} and (1) above,
\[
  \frac{d}{d \epsilon_r} \tilde{E}
 =\sum_{l=1}^{\infty} \sum_{i,j=1}^2 (r_l)_{ij} \mathfrak{r}_l^{ij} \tilde{E}
  =0,
\]
where $\mathfrak{r}_l^{ij}$ is the operation $\mathfrak{r}_l$, with
the two new half-edges called $i$ and $j$ (for the $i$-th and $j$-th points).
It remains to prove that $(r_l)_{ij}$ is ``non-degenerate'' in the sense that
$(r_l)_{ij} \ne 0$ unless 
\[
  \mathfrak{r}_l^{ij} \tilde{E} =0
\]
due to the \emph{a priori} constraint \eqref{e:r}.
The non-vanishing of the $r_l$-matrices follows from Proposition~\ref{p:A1}, 
which will be proved in the Appendix. 
Therefore, 
\[
  \mathfrak{r}_l^{ij} \tilde{E} = 0
\]
which implies \eqref{e:rl2}.
\end{proof}

\subsection{Some applications of Theorem~\ref{t:c1}} \label{s:6.5}

In a series of joint work with D.~Arcara \cite{AL1, AL2, AL3} and
with A.~Givental \cite{GL},
Theorem~\ref{t:c1} is shown to implies 
\begin{itemize}
\item A uniform derivation of all known $g=1, 2$ tautological equations.
\item Derivation of a new tautological equation in $\Mbar_{3,1}$ of 
codimension $3$.
\item All monomials of $\kappa$-classes and $\psi$-classes are 
independent in $R^k(\Mbar_{g,n})/R^k(\partial\Mbar_{g,n})$ 
for all $k \leq [g/3]$.
\end{itemize}

The first two investigate the \emph{existence} of tautological equations,
while the last one deals with the \emph{non-existence} of tautological 
equations.
Theorem~\ref{t:c1} provides a \emph{new and systematic} way to study
the existence or non-existence of tautological equations via
induction and linear algebra.
In fact, it is the simplicity of the derivations of the above results
that might point to some new direction in the study of tautological rings.

\appendix

\section{The $R$ matrix for $P^1$}
\begin{center} {\footnotesize by Y.~IWAO and Y.-P.~LEE} \end{center}

\vspace{10 pt}

The notations here follows those in \cite{LP}. 
All invariants and functions are for $X = \PP^1$.

\[
  R(z) = \sum_{n=0}^{\infty} R_n z^n
\]
is defined by $R_0 = I$ (the $2 \times 2$ identity matrix) and
the following recursive relation (1.4.5 in \cite{LP}), 
 \begin{multline}  \label{e:22}
 -\frac{\sqrt{-1}}{2}(u_1-u_2)^{-1}
  \left(\begin{array}{ccc}
  (R_{n-1})_{1}^{2} & (R_{n-1})_{2}^{2}\\
  (R_{n-1})_{1}^{1} & (R_{n-1})_{2}^{1}
       \end{array}\right)(du_1-du_2)\\\\
 +\left(\begin{array}{ccc}
  (dR_{n-1})_{1}^{1} & (dR_{n-1})_{2}^{1}\\
  (dR_{n-1})_{1}^{2} & (dR_{n-1})_{2}^{2}
       \end{array}\right)
 =\left(\begin{array}{ccc}
  0 & (R_{n})_{2}^{1}\\
  (R_{n})_{1}^{2} & 0
       \end{array}\right)(du_1-du_2) 
 \end{multline}
for $n \geq 1$, where $(u_1,u_2)$ are the canonical coordinates of 
$QH^*(\PP^1)$. 

\begin{proposition} \label{p:1}
 \begin{equation}\label{eq:R_n} 
  R_n=
  \left(\begin{array}{ccc}
   1 & 2n(-1)^{n-1}\sqrt{-1}\\
   2n\sqrt{-1} & (-1)^n
       \end{array}\right)\frac{c_n}{v^n},
 \end{equation}
where
 \begin{displaymath} 
  c_n:= -\frac{\prod_{k=1}^{n-1} (-1+4k^2)}{2^{2n}n!}, \qquad v:=u_1-u_2.
 \end{displaymath}
and it is understood that $c_1 = -1/4$.
\end{proposition}

\begin{proof}
Equation \eqref{e:22} gives a recursive relation which determines all
$R_n$ from $R_0 = I$.
It is easy to check that \eqref{eq:R_n} satisfies \eqref{e:22}. 
\end{proof}

Recall that $r(z)$ is defined as $\log R(z)$.

\begin{proposition} \label{p:A1}
  \begin{displaymath}
  r_l=
  \begin{cases}
   \left(\begin{array}{ccc}
    0 & -c_l\sqrt{-1}\\
    c_l\sqrt{-1} & 0 
       \end{array}\right)& \text{if $l$ is even},\\\\
   \left(\begin{array}{ccc}
    a_l & b_l\sqrt{-1}\\
    b_l\sqrt{-1} & -a_l 
       \end{array}\right)& \text{if $l$ is odd},
  \end{cases}
  \end{displaymath}
such that $a_l, b_l, c_l$ are all nonzero rational numbers.
\end{proposition}

The rest of the Appendix is used to prove this proposition.

First, it is not very difficult to see the matrices $r_l$ should be of
the above forms and all $a_l, b_l, c_l$ are rational numbers. 
These assertions follow from the following formula for $r_l$
\begin{equation} \label{e:24}
 \begin{split}
  \sum_{l=1}^{\infty} r_l z^l &= \log R(z) 
  = -\sum_{n=0}^{\infty} \dfrac{(1-R(z))^{n+1}}{n+1} \\
  &= \sum_l z^l \sum_{m=1}^l \dfrac{(-1)^{m-1}}{m} 
  \sum_{\substack{i_1+\cdots+i_m=l \\ i_j>0}} R_{i_1}\cdots R_{i_m},
 \end{split}
\end{equation}
Equation \eqref{e:r}, and induction on $l$.

Here are a few examples of $r_l$ matrices obtained from \eqref{e:24} and
Proposition~\ref{p:1}:
 \begin{align*}
  r_1& =
  \dfrac{1}{4}
  \left(\begin{array}{ccc}
    -1 & -2\sqrt{-1}\\
    -2\sqrt{-1} & 1 
       \end{array}\right)\\
  r_2& =
  \dfrac{1}{2^3}
  \left(\begin{array}{ccc}
    0 & 3\sqrt{-1}\\
    -3\sqrt{-1} & 0 
       \end{array}\right)\\
  r_3& =
  \dfrac{1}{2^5}
  \left(\begin{array}{ccc}
    -4 & -23\sqrt{-1}\\
    -23\sqrt{-1} & 4 
       \end{array}\right)\\ 
  r_4& =
  \dfrac{1}{2^4}
  \left(\begin{array}{ccc}
    0 & 33\sqrt{-1}\\
    -33\sqrt{-1} & 0 
       \end{array}\right)\\
  r_5& =
  \dfrac{1}{5\cdot 2^9}
  \left(\begin{array}{ccc}
    -2132 & -20839\sqrt{-1}\\
    -20839\sqrt{-1} & 2132 
       \end{array}\right).
 \end{align*}

The non-vanishing of $a_l, b_l, c_l$ requires some elementary 
(though somehow lengthy) arguments.
The idea is to get some rough estimates of the absolute values of entries of 
$r_l$, which are enough to guarantee the non-vanishing.
We will start with some preparations.

Define a function of $n$, $C(n)$, by
 \begin{displaymath}
  C(n)= 
  \prod_{k=1}^{n-1} (1-\dfrac{1}{4k^2}). 
  \end{displaymath}
Then for $n \geq 1$,
 \begin{displaymath}
  (R_n)_{1}^{2}=C(n)\cdot (n-1)!\cdot (\dfrac{-\sqrt{-1}}{2})\cdot v^{-n}.
 \end{displaymath}

\begin{lemma} \label{l:5}
$C(n)$ is strictly decreasing and $\lim_{n \rightarrow \infty} C(n) > 0.62$.
\end{lemma}

\begin{proof}
$C(n)$ is obviously strictly decreasing,
so does $\lim_{n \rightarrow \infty} [\ln C(n)]= \ln [\lim_{n \rightarrow \infty} C(n)]$.
Since 
 \begin{displaymath}
  \ln (1-x) \geq 4(\ln \dfrac{3}{4})x \qquad \textrm{for }0< x \leq \dfrac{1}{4}, 
 \end{displaymath}
\[
  \lim_{n \rightarrow \infty} [\ln C(n)]= 
\lim_{n \rightarrow \infty} \left[\sum_{k=1}^{n-1} \ln(1-\dfrac{1}{4k^2}) \right] 
\geq \dfrac{\pi^2}{6} \ln \dfrac{3}{4}.
\]
Therefore
 \begin{displaymath}
  \lim_{n \rightarrow \infty} C(n) \geq e^{\frac{\pi^2}{6} \ln \frac{3}{4}}
  \approx 0.62299 \cdots.
 \end{displaymath}
\end{proof}

By abuse of notation, 
\emph{let $|(R_n)_{i}^{j}|$ denote the absolute value of the coefficient of 
$v^{-n}$ in $(R_n)_{i}^{j}$}.
Since $\displaystyle{(R_n)_{1}^{1}= \frac{-\sqrt{-1}}{2n}(R_n)_{1}^{2}}$,  we get the following corollary:

\begin{corollary} \label{c:3}
  For $n \geq 1$,
  \begin{align*}
   &0.62 \dfrac{(n-1)!}{2} < |(R_n)_{1}^{2}| \leq \dfrac{(n-1)!}{2}\\
   &0.62 \dfrac{(n-1)!}{4n} < |(R_n)_{1}^{1}| \leq \dfrac{(n-1)!}{4n}. 
  \end{align*}
\end{corollary}

Now define $\sigma_{n}^{l}$ by
 \begin{displaymath}
  \sigma_{n}^{l}:= \sum_{\substack{i_1+\cdots +i_l=n\\i_j>0}} i_1!i_2!\cdots i_l!,
 \end{displaymath} 

\begin{lemma} \label{c:4}
  \begin{displaymath}
   \sigma_{n}^{l} \leq (\frac{8}{3})^{l-1}n! \qquad \text{for } n=0, 1, 2, 3,\cdots, \quad l=1, 2, 3,\cdots.
  \end{displaymath}
\end{lemma}

\begin{proof}
It is obviously true for $l=1$ and arbitrary $n$.
By induction on $n$ with fixed $l=2$:
  \begin{displaymath}
   \sigma_{n}^{2} \leq \dfrac{8}{3}n! \qquad n=0, 1, 2, 3,\cdots.
  \end{displaymath}
Then induction on $l$.
Assume $\sigma_{n}^{l} \leq (\dfrac{8}{3})^{l-1}n!$ for all $n$.  
It suffices to show $\sigma_{n}^{l+1} \leq (\dfrac{8}{3})^{l}n!$ for all $n$.
 \begin{align*}
  \sigma_{n}^{l+1}& =0!\sigma_{n}^{l}+1!\sigma_{n-1}^{l}+2!\sigma_{n-2}^{l}+\cdots +(n-1)!\sigma_{1}^{l}+n!\sigma_{0}^{l}\\
& \leq (\dfrac{8}{3})^{l-1}[0!n!+1!(n-1)!+2!(n-2)!+\cdots + (n-1)!1!+n!0!]\\
& =(\frac{8}{3})^{l-1}\sigma_{n}^{2} \leq (\dfrac{8}{3})^{l}n!.
 \end{align*} 
\end{proof}

\begin{lemma} \label{l:9}
  Let 
  \begin{displaymath}
   S(m):=\sum_{k=3}^{m}\dfrac{2^{k-3}(m-k)!}{k} \quad \textrm{for $m \geq 5$}.
  \end{displaymath}

Then
\begin{displaymath}
    T(m):=\dfrac{S(m)}{(m-3)!} \quad \textrm{is a decreasing function of $m$ for $m \geq 5$}.
\end{displaymath}
\end{lemma}

\begin{proof}
Clearly, $i$th term of $T(m) \geq$ $i$th term of $T(m+1)$ for $i=1, 2, \cdots m-3$, with equality only for $i=1$. 
For $m \geq 5$, $T(m+1) < T(m)$ as
\[
 \dfrac{2^{m-3}}{m(m-2)!}+ \dfrac{2^{m-2}}{(m+1)(m-2)!} < \dfrac{2^{m-3}}{m(m-3)!}.
\]
\end{proof}

\begin{corollary} \label{c:5}
  \begin{enumerate}
   \item $T(m) < 1$ \quad {\it for $m \geq 5$}.
   \item $T(m) < 0.477$ \quad {\it for $m \geq 9$}.
   \item $S(m) < (m-3)!$ \quad {\it for $m \geq 5$}.
   \item $S(m) < 0.477(m-3)!$ \quad {\it for $m \geq 9$}.
  \end{enumerate}
\end{corollary} 

\begin{corollary} \label{c:6}
  \begin{displaymath}
   \dfrac{1}{2} \left| \sum_{k=1}^{m-1} (R_kR_{m-k})_{1}^{1} \right| \leq \dfrac{15}{256}(m-3)! \quad \textrm{for $m$ odd, $m \geq 7$}.
  \end{displaymath}
\end{corollary}

\begin{proof}
Since $m$ is odd, 

 \begin{align*}
  &\dfrac{1}{2} \left| \sum_{k=1}^{m-1} (R_kR_{m-k})_{1}^{1} \right| 
 \leq \sum_{k=1}^{\frac{m-1}{2}} \left|(R_k)_{1}^{1} \right| \left|R_{m-k})_{1}^{1} \right|\\
  \leq &\dfrac{1}{4}\cdot \dfrac{3(m-2)!}{16(m-1)}+\dfrac{3\cdot 1!}{16 \cdot2}\cdot \dfrac{3(m-3)!}{16(m-2)} 
+\cdots +\dfrac{3(\dfrac{m-3}{2})!}{16(\dfrac{m-1}{2})}\cdot \dfrac{3(\dfrac{(m-1)}{2})!}{16(\dfrac{(m+1)}{2})}\\
  \leq &\dfrac{3}{64}(m-3)!+ \dfrac{9}{256}[0!(m-4)!+1!(m-5)!+\cdots +(\dfrac{m-5}{2})!(\dfrac{m-3}{2})!\\
  \leq &\dfrac{3}{64}\cdot \dfrac{5}{4}\cdot (m-3)!= \dfrac{15}{256}(m-3)! \quad \textrm{for $m \geq 7$}.
  \end{align*}
\end{proof}

\begin{lemma} \label{l:11}
  \begin{displaymath}
   \left|(R_m)_{1}^{1} \right| \geq 0.13(m-2)! \quad \textrm{for $m \geq 7$}.
  \end{displaymath}
\end{lemma}

\begin{proof}
By Corollary~\ref{c:3}

 \begin{align*}
  \left|(R_m)_{1}^{1} \right| &> 0.62\cdot \dfrac{(m-1)!}{4m}=\frac{1}{4}\cdot 0.62\cdot (1-\dfrac{1}{m})(m-2)!\\
  &\geq \frac{1}{4}\cdot 0.62\cdot \dfrac{6}{7}(m-2)!\\
  &> 0.13(m-2)!.
 \end{align*} 

\end{proof}

Now we are ready to prove Proposition~\ref{p:A1}. 
We start with off-diagonal entries.

\begin{displaymath}
  (r_l)_{1}^{2}= (R_l)_{1}^{2}+\underbrace{\sum_{n=2}^l \dfrac{(-1)^{n-1}}{n} \sum_{\substack{i_1+\cdots+i_n=l \\ i_j>0}} (R_{i_1}\cdots R_{i_n})_{1}^{2}}_{=:(R'_l)_{1}^{2}}.
\end{displaymath}

By the triangular inequality, we have
\begin{displaymath}
  |(r_l)_{1}^{2}| \geq \left| \left|(R_l)_{1}^{2} \right|- \left|(R'_l)_{1}^{2} \right| \right|.
\end{displaymath}
So it suffices to show
 \begin{displaymath}
  \left|(R_l)_{1}^{2} \right| > \left|(R'_l)_{1}^{2} \right|
 \end{displaymath}

Now,

\begin{align*}
  &\left|(R'_l)_{1}^{2} \right| \\
\leq &\sum_{n=2}^l \dfrac{1}{n} \sum_{\substack{i_1+\cdots+i_n=l \\ i_j>0}} \left| (R_{i_1}\cdots R_{i_n})_{1}^{2} \right| \\
  \leq &\sum_{n=2}^l \dfrac{1}{n} \sum \dfrac{(i_{j_1}-1)!}{2}\cdots \dfrac{(i_{j_k}-1)!}{2}\cdot \dfrac{(i_{j_{k+1}}-1)!}{4i_{j_{k+1}}}\cdots \dfrac{(i_{j_n}-1)!}{4i_{j_n}} \quad \textrm{(by Corollary~\ref{c:3})}\\
 \leq &\sum_{n=2}^l \dfrac{1}{n}\cdot \dfrac{3^n-(-1)^n}{2\cdot 4^n}\cdot (\dfrac{8}{3})^{n-1}(l-n)! \quad \textrm{(by Lemma~\ref{c:4})}\\
  = &\dfrac{3}{16} \sum_{n=2}^l \dfrac{1}{n} \left[ 2^n-(-\dfrac{2}{3})^n \right](l-n)! \\
  < &\dfrac{3}{16} \sum_{n=2}^l \dfrac{1}{n}\cdot 2^n\cdot (l-n)! \\
  < &\dfrac{3}{8}(l-2)!+ \dfrac{3}{2}(l-3)!  \quad \textrm{for $l \geq 1$} \quad \textrm{(by Corollary~\ref{c:5})}.
 \end{align*}

By Corollary~\ref{c:3},
\[
 0.31(l-1)!< \left| (R_l)_{1}^{2}\right| < 0.5(l-1)! \quad \text{for $l \geq 1$.}
\]
So for $l \geq 5$,
 \begin{displaymath}
  \left| (R_l)_{1}^{2}\right|> 0.3(l-1)!> \dfrac{3}{8}(l-2)!+ \dfrac{3}{2}(l-3)!> \left| (R'_l)_{1}^{2}\right|,
 \end{displaymath}
i.e., $(r_l)_{1}^{2} \neq 0$.
Since we know $(r_l)_{1}^{2} \neq 0$ for $l=1, 2, 3, 4$, 
$(r_l)_{1}^{2} \neq 0$ for all $l \geq 1$.

Similarly for diagonal terms for $l=odd$. 

 \begin{align*}
  \left| (R'_l)_{1}^{1}\right| 
  &\leq \dfrac{1}{2} \sum_{\substack{i_1+i_2=l\\i_j>0}} \left| (R_{i_1}R_{i_2})_{1}^{1}\right|+ \sum_{n=3}^l \dfrac{1}{n} \sum_{\substack{i_1+\cdots+i_n=l \\ i_j>0}} \left| (R_{i_1}\cdots R_{i_n})_{1}^{1} \right|\\
  &\leq \dfrac{15}{256}(l-3)!+ \dfrac{3}{16} \sum_{n=3}^l \dfrac{1}{n} \left[ 2^n+(-\dfrac{2}{3})^n \right](l-n)!\\
  &< \dfrac{15}{256}(l-3)!+ \dfrac{3}{16} \sum_{n=3}^l \dfrac{1}{n}\cdot 2^n\cdot (l-n)!\\
  &= \dfrac{15}{256}(l-3)!+ \dfrac{3}{2} S(l)\\
  &< \dfrac{15}{256}(l-3)!+ 0.72(l-3)! \quad \textrm{for $l \geq 9$, $m=$ odd} \quad \textrm{(by Corollary~\ref{c:5})}.
 \end{align*}
By Lemma~\ref{l:11},
 \begin{displaymath}
  \left| (R_l)_{1}^{1}\right|> 0.13(l-1)!> (\dfrac{15}{256}+ 0.72)(l-3)!> \left| (R'_l)_{1}^{1}\right| \quad \textrm{for $l \geq 9$}.
 \end{displaymath}
$(r_l)_{1}^{1} \neq 0$ for $l \geq 9$, $l=$ odd.
We know $(r_l)_{1}^{1} \neq 0$ for $l=1, 3, 5$.

It is now left to check $l=7$, which is checked by hand.
 \begin{align*}   
  \left| (R'_7)_{1}^{1}\right| &\leq \dfrac{15}{256}\cdot 4!+ \sum_{m=3}^{7}\dfrac{3^m+(-1)^m}{2\cdot 4^m}\cdot \sigma_{7-m}^{m}\\
  &= \dfrac{15\cdot 24}{256}+ \dfrac{1}{2} \left[ \dfrac{1}{3}\cdot \dfrac{27-1}{64}\cdot 108+ \dfrac{1}{4}\cdot \dfrac{81+1}{256}\cdot 52 \right.\\
  &+ \left. \dfrac{1}{5}\cdot \dfrac{243-1}{1024}\cdot 20+ \dfrac{1}{6}\cdot \dfrac{729+1}{4096}\cdot 6+ \dfrac{1}{7}\cdot \dfrac{2187-1}{16384}\cdot 1 \right]= 11.3720\cdots.
 \end{align*}
On the other hand,
\begin{displaymath}
  \left| (R_7)_{1}^{1}\right|= \dfrac{3\cdot 15\cdot 35\cdot 63\cdot 99\cdot 143}{2^{14}\cdot 7!}= 17.0114\cdots.
 \end{displaymath}
So $(r_7)_{1}^{1} \neq 0$.
The proof of Proposition~\ref{p:A1} is now complete.

\end{document}